\input amstex
\input epsf
\documentstyle{amsppt}
\loadeurm
\hfuzz = 1,1pt

\def\R{\Bbb R}

\def\C{\Bbb C}
\def\H{\Bbb H}
\def\A{\Cal A}
\def\ZZ{\Cal Z}

\def\X{\Cal X}

\def\TE{\Cal T}
\def\pTE{\widetilde{\Cal T}}

\def\ls{\leqslant}
\def\gs{\geqslant}

\def\diag{\operatorname{diag}}

\def\alt{\operatorname{alt}}
\def\spa{\operatorname{span}}
\def\corank{\operatorname{corank}}

\def\Edge{\operatorname{Edge}}
\def\Vert{\operatorname{Vert}}

\def\Im{\operatorname{Im}}

\topmatter

\title{Cluster algebras and Weil-Petersson forms}
\endtitle
\author Michael Gekhtman$^*$, Michael Shapiro$^\ddag$,
 and Alek Vainshtein$^\dag$ \endauthor
\affil $^*$ Department of Mathematics, University of Notre Dame,
Notre Dame, IN 46556,
{\tt Michael.Gekhtman.1\@nd.edu}\\
$^\ddag$ Department of Mathematics, Michigan State University,
East Lansing, MI 48824-1027, {\tt mshapiro\@math.msu.edu}\\
$^\dag$ Departments of Mathematics and of Computer Science,
University of Haifa,  Israel 31905,
{\tt alek\@math.haifa.ac.il}
\endaffil
\leftheadtext{Michael Gekhtman, Michael Shapiro, and Alek Vainshtein}
\endsubjclass
\thanks{This work was started during the visit of the
second and the third authors to Max-Planck Institute f\"ur Mathematik in 
Bonn in summer 2002, and finished during their visit to MSRI, Berkeley, 
in spring 2004. 
The authors are partially supported by the BSF grant 2002375.
The second author is partially supported by the Institute of Quantum Sciences,
Michigan State University.}
\endthanks

\abstract In the previous paper~\cite{GSV} we have discussed
Poisson properties of cluster algebras of geometric type for the
case of a nondegenerate matrix of transition exponents. In this
paper we consider the case of a general matrix of transition
exponents. Our leading idea is that a relevant geometric object in
this case is a certain closed 2-form  compatible with the cluster
algebra structure. The main example is provided by Penner
coordinates on the decorated Teichm\"uller space, in which case
the above form coincides with the classical Weil-Petersson
symplectic form.
\endabstract

\endtopmatter

\heading{1. Introduction}
\endheading

Cluster algebras are a class of commutative rings defined
axiomatically in terms of a distinguished family of generators
called cluster variables. The study of cluster algebras was
started by S.~Fomin and A.~Zelevinsky in~\cite{FZ2} and then
continued in~\cite{FZ3},~\cite{FZ4},~\cite{GSV}. The main
motivation for cluster algebras came from the study of dual
canonical bases, see~\cite{BZ},~\cite{Ze2} and the theory of
double Bruhat cells in complex semisimple Lie groups,
see~\cite{SSV1},~\cite{SSV2},
~\cite{FZ1},~\cite{SSVZ},~\cite{Ze1}.

This paper is a continuation of~\cite{GSV}, where we have
developed a Poisson formalism compatible with cluster algebra
structures. Here we suggest to associate with a cluster algebra a
compatible skew-symmetric $2$-form. Under certain irreducibility
assumptions this $2$-form is unique up to a scalar factor. In a
non-degenerate case it is a symplectic form dual to the unique
compatible Poisson structure constructed in~\cite{GSV}. In the
general case there exists a projection of the cluster manifold
defined in~\cite{GSV} to a rational manifold of a smaller
dimension equipped with special $\tau$-coordinates and their
cluster transformations. On the latter manifold one has a natural
compatible Poisson structure and a symplectic form dual to each
other. Applying this construction to the coordinate ring of the
Teichm\"uller space we recover the famous Weil-Petersson
symplectic form. The initial cluster in this case is formed by
{\it Penner coordinates} \cite{Pe} associated with an ideal
triangulation of a punctured Riemann surface. These coordinates
proved to be a crucial  ingredient in the study of decorated
Teichm\" uller spaces \cite{Pe}, \cite{Fo} and their quantization
\cite{Ka}. Motivated by this important example, we call the
compatible symplectic form on the rational manifold associated
with any cluster algebra {\it the Weil-Petersson symplectic form}.

Since the earlier version of this paper was posted, two preprints
by  Fock and Goncharov \cite{FG1}, \cite{FG2} has appeared that explored
cluster algebra aspects of the moduli space of local systems on a
punctured Riemann surface.

The authors would like to thank V.~Fock, S.~Fomin, A.~Goncharov,
A.~Postnikov, and A.~Zele\-vinsky for useful discussions,
N.~Ivanov for crucial consultations on triangulations of surfaces,
and D.~Thurston for explaining to us his unpublished result.
Sincere thanks go to anonymous referees for a thorough reading of
the manuscript and advice on how to improve it;
in particular, for pointing out gaps in the original proof
of Theorem 3.4 and for suggesting a shorter proof of Theorem 3.3.

\heading{2. Degenerate cluster algebras of geometric type}
\endheading

Let $\A$ be a cluster algebra of geometric type. As it has been
discussed in~\cite{FZ2}, $\A$ is completely determined by the
initial cluster with a fixed set of cluster variables
$x_1,\dots,x_n$ and a skew-symmetrizable transformation matrix $Z$. In
our previous paper~\cite{GSV} we have considered the case of
nondegenerate $Z$. This paper discusses the case of general $Z$.
For the sake of simplicity we consider only the skew-symmetric
case.

We start by recalling main definitions
from~\cite{FZ2} and~\cite{GSV}.

\subheading{2.1. Cluster algebras of meromorphic functions on a rational
manifold}
Let $\ZZ_n$ be the set of integer-valued skew-symmetric $n\times n$
matrices. According to~\cite{FZ2}, any $Z=(z_{ij})\in\ZZ_{n}$ defines a
cluster algebra of geometric type in the following way. Let us fix
a set of $n$ {\it cluster\/} variables $f_1,\dots,f_n$. For each
$i\in [1,n]$ we introduce a transformation $T_i$ of cluster
variables by
$$
T_i(f_j)=\bar
f_j=\left\{\alignedat2 &\frac1{f_i}\left(\prod_{z_{ik}>0}
f_k^{z_{ik}}+\prod_{z_{ik}<0}
  f_k^{-z_{ik}}\right)\quad& &\text{for $j= i$}\\
&f_j \qquad& &\text{for $j\ne i$},
\endalignedat\right. \tag 2.1
$$
and the corresponding matrix transformation $\bar Z=T_i(Z)$,
called {\it mutation}, by
$$
\bar z_{kl}=\left\{\alignedat2
&-z_{kl}\qquad& &\text{for $(k-i)(l-i)=0$}\\
&z_{kl}+\frac{|z_{ki}|z_{il}+z_{ki}|z_{il}|}2\quad& &\text{for
$(k-i)(l-i)\ne
  0$}.
\endalignedat\right. \tag 2.2
$$
Observe that $T_i$ is an involution and
$\bar Z$ belongs to $\ZZ_{n}$. Thus, one can apply
transformations $T_i$ to the new set of cluster variables (using
the new matrix), etc. The {\it rank $n$ cluster algebra\/} 
({\it of geometric type\/})
is the subalgebra of the field of rational functions in cluster
variables $f_1,\dots,f_n$ generated by the union of all clusters;
its ground ring is the ring of integers.
We denote this algebra by $\A(Z)$.

One can represent $\A(Z)$ with the help of an $n$-regular tree
$\Bbb T_n$ whose edges are labeled by the numbers $1,\dots,n$ so
that the $n$ edges incident to each vertex receive different
labels. To each vertex $v$ of $\Bbb T_n$ we assign a set of $n$
variables $f_{v,1},\dots,f_{v,n}$, which together form 
a {\it cluster\/} at $v$. For an edge $v\bar
v$ of $\Bbb T_n$ that is labeled by $i\in [1,n]$, the clusters
$f=f_v$ and $\bar f=f_{\bar v}$ are related by the transformation
$T_i$ given by~(2.1).

\example{Example 1}  
Consider a $3\times 3$ skew-symmetric matrix
$$
Z=\pmatrix \phantom{-} 0 & \phantom{-}1 & \hfill-1 \cr
           \hfill-1& \phantom{-}0 &  \phantom{-}1 \cr
           \phantom{-}1 & \hfill-1&  \phantom{-}0 \endpmatrix
$$
and fix the initial cluster $f_1, f_2, f_3$. Then the  
action of transformations $T_1, T_2, T_3$ is given by
$$\alignat4
T_1(f_1) &= \frac{f_2+f_3}{f_1},&\quad T_1(f_2)&=f_2,& \quad T_1(f_3)&=f_3,&
\qquad T_1(Z)&=
\pmatrix \phantom{-}0 & \hfill-1 & \phantom{-}1 \cr
               \phantom{-}1& \phantom{-}0 &  \phantom{-}0 \cr
               \hfill-1 & \phantom{-}0&  \phantom{-}0 \endpmatrix;\\
T_2(f_2) &= \frac{f_1+f_3}{f_2},& \quad T_2(f_1)&=f_1, &\quad
T_2(f_3)&=f_3,&\qquad T_2(Z)&=
\pmatrix \phantom{-}0 & \hfill-1 & \phantom{-}0 \cr
               \phantom{-}1& \phantom{-}0 &  \hfill-1 \cr
               \phantom{-}0 & \phantom{-}1&  \phantom{-}0 \endpmatrix;\\ 
T_3(f_3) &= \frac{f_1+f_2}{f_3},&\quad T_3(f_1)&=f_1;&\quad T_3(f_2)&=f_2,
&\qquad T_3(Z)&=
\pmatrix \phantom{-}0 & \phantom{-}0 & \phantom{-}1 \cr
               \phantom{-}0& \phantom{-}0 &  \hfill-1 \cr
               \hfill-1 & \phantom{-}1&  \phantom{-}0 \endpmatrix.
\endalignat
$$
\endexample

Consider an $n$-dimensional rational manifold $\Cal M$ (recall that this
means that the field $F(\Cal M)$ of meromorphic functions on $\Cal M$ is a
transcendental extension of the ground field, i.e. the field of rational
functions in $n$ variables).
Let us fix an arbitrary set of $n$ generators $\{x_1,\dots,x_n\}$ of the
field $F(\Cal M)$. Then we get an
embedding of the cluster algebra $\A(Z)$ into $F(\Cal M)$ by mapping
the initial cluster variables $f_1,\dots, f_n$ to $x_1,\dots, x_n$,
respectively. In what follows we assume that this embedding is fixed and
identify the elements of $\A(Z)$ with the corresponding meromorphic
functions on $\Cal M$.

Recall that in our previous paper~\cite{GSV} we have considered a
slightly different situation: only a part of the variables is
regarded as cluster variables; the transformations of the
remaining variables are prohibited, and these variables are called
{\it stable}. In this situation only $m$ rows of the
transformation matrix are given, where $m$ is the number of
cluster variables. An important nondegeneracy condition assumed
in~\cite{GSV} is that the rank of the corresponding $m\times n$
matrix  equals $m$ (it is proved in~\cite{GSV} that the rank is
invariant under mutations).

In our present situation, when $Z$ is an $n\times n$ matrix of
rank $m\ls n$, one can choose  $m$ rows of $Z$ so that the
resulting submatrix of $Z$ is nondegenerate, and  obtain the
situation described in~\cite{GSV}. The corresponding cluster
algebra is then a subalgebra of $\A(Z)$; it is denoted $\A_I(Z)$, where
$I=(i_1,\dots,i_m)$ is the multiindex corresponding to the chosen rows.

On the other hand, starting from a cluster algebra with both
cluster and stable variables, one can get a coefficient-free
cluster algebra (the one without stable variables) by just
fixing the values of all the stable variables at $1$.

\subheading{2.2. Compatible Poisson brackets}
Let $\{\cdot,\cdot\}$ be a Poisson bracket on an $n$-dimensional
algebraic manifold $\Cal M$.
We say that functions $g_1,\dots,g_n$ are {\it log-canonical}
with respect to $\{\cdot,\cdot\}$ if $\{g_i,g_j\}=p_{ij} g_ig_j$,
where $p_{ij}$ are constants. The matrix
$(p_{ij})$ is called the {\it coefficient matrix\/} of
$\{\cdot,\cdot\}$ (in the basis $g$); evidently, it is skew-symmetric.
Assume that $\Cal M$ is a rational manifold.
We say that a Poisson bracket $\{\cdot,\cdot\}$ on $\Cal M$
is {\it compatible\/}
with a cluster algebra $\A(Z)$ if all clusters in $\A(Z)$ are
log-canonical with respect to $\{\cdot,\cdot\}$.

In~\cite{GSV} we studied the space of Poisson brackets compatible
with algebras of type $\A_I(Z)$. In particular, we have proved
that such brackets (considered up to a scalar factor) form a
vector space of dimension $1+\binom{n-m}2$, where $m$ is the rank
of $Z$.
 As a consequence,
when $Z$ is nondegenerate, there exists a unique (up to a scalar factor)
Poisson bracket  compatible with $\A(Z)$.

Examples show that when a cluster algebra is associated with a
geometric object, the corresponding compatible Poisson brackets
also have a natural invariant interpretation. This is the case,
e.g., for cluster algebras associated with double Bruhat cells in
a semi-simple Lie group. There, the initial cluster can be chosen
to consist of certain generalized minors. These minors form a
log-canonical coordinate system with respect to the standard
Poisson-Lie structure on the group. Moreover, the corresponding
Poisson bracket is compatible with cluster transformations
(see~\cite{FZ1},~\cite{KZ} for definitions and details). Another
example, discussed in~\cite{GSV}, is the real Grassmannian
$G(k,n)$, viewed as a Poisson homogeneous space of the group
$SL(n)$ equipped with the Sklyanin Poisson-Lie bracket. A
log-canonical coordinate system for a push-forward of the Sklyanin
bracket to  $G(k,n)$ is formed by a subset of Pl\"ucker
coordinates. These coordinates serve as an initial cluster for the
cluster algebra compatible with this Poisson bracket. It is proved
in~\cite{GSV} that the resulting cluster algebra contains all  the
Pl\"ucker coordinates and thus, the coordinate ring of $G(k,n)$ is
endowed with a cluster algebra structure.

However, the following example shows that a bracket compatible
with $\A(Z)$ may not exist if $Z$ is degenerate.

\example{Example 2} 
Let $\{\cdot,\cdot\}$ be a Poisson structure compatible with the
cluster algebra described above in Example~1. 
Then, by the definition of compatibility, 
$\{f_1,f_2\}=\lambda f_1f_2$, $\{f_2,f_3\}=\mu
f_2f_3$, $\{f_1,f_3\}=\nu f_1f_3$. Applying $T_1$ one gets 
$\{\bar f_1,f_2\}=-\lambda f_2^2/f_1-\mu
f_2f_3/f_1 -\lambda f_2f_3/f_1$. On the other hand, the compatibility yields
$\{\bar f_1,f_2\}=\alpha \bar f_1f_2=\alpha f_2^2/f_1+\alpha f_2f_3/f_1$. 
Comparing
these two expressions we immediately get $\mu=0$, which means that
$\{f_2,f_3\}=0$. Similarly, $\{f_1,f_2\}=\{f_1,f_3\}=0$. Hence we see that the
only Poisson structure compatible with this cluster algebra is
trivial.
\endexample

\subheading{2.3. Compatible 2-forms}
Motivated by the example above, we will look for an alternative
geometric object compatible with a general cluster algebra. Natural
duality suggests a compatible differential 2-form as a possible candidate.

Let $\omega$ be a differential 2-form on an $n$-dimensional algebraic manifold.
We say that functions $g_1,\dots,g_n$ are {\it log-canonical}
with respect to $\omega$ if
$$
\omega=\sum_{i,j=1}^n\omega_{ij}\frac{d g_i\wedge d g_j}{g_i g_j},
$$
where $\omega_{ij}$ are
constants. The matrix $\Omega=(\omega_{ij})$ is called the {\it
coefficient matrix\/} of $\omega$ (in the basis $g$); evidently,
$\Omega$ is skew-symmetric. We say that a 2-form $\omega$ on a rational
manifold is {\it
compatible\/} with a cluster algebra $\A(Z)$ if all clusters in
$\A(Z)$ are log-canonical with respect to $\omega$.

We say that a skew-symmetric matrix $M$ is {\it reducible\/} if there
exists a permutation matrix $P$ such that $PMP^T$ is a block-diagonal
matrix, and {\it irreducible\/} otherwise;
$r(M)$ is defined as the maximal number of diagonal blocks in $PMP^T$.
The partition into blocks defines an obvious equivalence relation
$\sim$ on the rows (or columns) of $M$.

\proclaim{Theorem 2.1}
The 2-forms on a rational $n$-dimensional manifold
compatible with $\A(Z)$ form a vector space of
dimension $r(Z)$.
Moreover, the coefficient matrices of these 2-forms
in the basis formed by the cluster variables are characterized by the
equation $\Omega=\Lambda Z$,
where $\Lambda=\diag(\lambda_1,\dots,\lambda_n)$ with
$\lambda_i=\lambda_j$ whenever $i\sim j$.
In particular, if $Z$ is irreducible, then
$\Omega=\lambda Z$.
\endproclaim

\demo{Proof}
Indeed, let $\omega$ be a 2-form compatible with $\A(Z)$. Then
$$
\omega=\sum_{j,k=1}^n\omega_{jk}\frac{df_j\wedge d f_k}{f_j f_k}
=\sum_{j,k=1}^n\bar\omega_{jk}
\frac{d\bar f_j\wedge d\bar f_k}{\bar f_j\bar f_k},
$$
where $\bar f_j=T_i(f_j)$ is given by~(2.1)
and $\bar\omega_{jk}$ are the coefficients of
$\omega$ in the basis $\bar f$. Recall that the only variable in the cluster
$\bar f$ that is different from the corresponding variable in $f$ is
$\bar f_i$, and
$$
\multline
d\bar f_i=-\frac1{f_i^2} \left(\prod_{z_{ik}>0}f_k^{z_{ik}}
+\prod_{z_{ik}<0}f_k^{-z_{ik}}\right)d f_i+\\
\sum_{z_{ik}>0}\frac{z_{ik}}{f_k}\dfrac1{1+\prod_{k=1}^nf_k^{-z_{ik}}}d f_k-
\sum_{z_{ik}<0}\frac{z_{ik}}{f_k}\dfrac1{1+\prod_{k=1}^nf_k^{z_{ik}}}d f_k.
\endmultline
$$
Thus, for any $j\ne i$ we immediately get
$$
\bar\omega_{ij}=-\omega_{ij}. \tag 2.3a
$$

Next, consider any pair $j,k\ne i$ and assume that both $z_{ij}$ and $z_{ik}$
are positive. Then
$$
\omega_{jk}=\bar\omega_{jk}+
\frac{\bar\omega_{ik}z_{ij}+\bar\omega_{ji}z_{ik}}
{1+\prod_{k=1}^nf_k^{-z_{ik}}}.
$$
This equality can only hold if $\bar\omega_{ik}z_{ij}+\bar\omega_{ji}z_{ik}=0$,
which gives
$$
\frac{\omega_{ij}}{z_{ij}}=\frac{\omega_{ik}}{z_{ik}}=\lambda_{i}.
\tag 2.4
$$
Besides, in this situation $$
\bar\omega_{jk}=\omega_{jk}. \tag 2.3b
$$

Similarly, if both $z_{ij}$ and $z_{ik}$ are negative, then
$$
\omega_{jk}=\bar\omega_{jk}-
\frac{\bar\omega_{ik}+\bar\omega_{ji}z_{ik}}
{1+\prod_{k=1}^nf_k^{z_{ik}}},
$$
and hence the same two relations as above are true.

Finally, let
$z_{ij}$ and $z_{ik}$ have different signs, say, $z_{ij}>0$ and
$z_{ik}<0$, then
$$
\omega_{jk}=\bar\omega_{jk}+
\frac{\bar\omega_{ik}z_{ij}\prod_{k=1}^nf_k^{z_{ik}}-\bar\omega_{ji}z_{ik}}
{1+\prod_{k=1}^nf_k^{z_{ik}}},
$$
which again leads to~(2.4); the only difference is that in this case
$$
\bar\omega_{jk}=\omega_{jk}+\omega_{ik}z_{ij}. \tag 2.3c
$$

By~(2.4), $\Omega=\Lambda Z$, where $\Lambda=\diag(\lambda_1,\dots,\lambda_n)$.
Since both $\Omega$ and $Z$ are skew-symmetric, we immediately get
$\lambda_i=\lambda_j$ whenever $i\sim j$.

It is worth mentioning that relations~(2.3a-c) are equivalent
to~(2.2). \qed
\enddemo

\subheading{2.4. $\tau$-coordinates and the Weil-Petersson form} 
Assume now that $\Cal M=\Cal X$ is the cluster manifold defined in~\cite{GSV}.
Recall that for a rank $n$ cluster algebra over a field $k=\R$ or $\C$, 
 $\Cal X$ is obtained by gluing together toric charts $(k^*)^n$ corresponding
to clusters. The corresponding cluster variables give the distinguished 
coordinate system in each chart; the gluing is defined by (2.1).
By~\cite{GSV, Lemma~2.1}, $\Cal X$ is rational.
Let $\omega$ be one of the forms described in Theorem~2.1. Note that $\omega$
is degenerate whenever $m<n$. Kernels of $\omega$ form an integrable
distribution (since it becomes linear in the coordinates $\log f$), so one
can define the quotient mapping $\pi^\omega\:\Cal X\to \Cal X^\omega$, 
where $\Cal X^\omega$
is a rational manifold of dimension $m$. Observe that this
mapping in fact does not depend on the choice of $\omega$, provided $\omega$ is
generic (that is, the corresponding diagonal matrix $\Lambda$ is nonsingular),
since the kernels of all generic $\omega$ coincide. The symplectic reduction
$\pi^\omega_* \omega$ of $\omega$ is thus well defined; moreover, 
$\pi^\omega_*\omega$
is a symplectic form on the image $\Cal X^\omega$.

Let us introduce coordinates in the image that are naturally related to
the initial coordinates $f$ on $\Cal X$. To do this, we define
another $n$-tuple of rational functions on $\Cal X$; it is denoted
 $\tau=(\tau_1,\dots,\tau_n)$ and is given by
$$
\tau_j=\prod_{k=1}^nf_k^{z_{jk}},\qquad j\in [1,n].
$$
Given a cluster $f$, we say that the entries of the corresponding
$\tau$ form a $\tau$-cluster.
It is easy to see that elements of the $\tau$-cluster are no longer
functionally independent.
To get a functionally independent
subset one has to choose at most $m$ entries in such a way that
the corresponding rows of $Z$ are linearly independent. Transformations
of $\tau$-coordinates corresponding to cluster transformations (2.1)
are studied in~\cite{GSV}.

Clearly,  any choice of $m$ functionally independent coordinates $\tau$
provides log-canonical coordinates for the 2-form $\pi^\omega_*\omega$.
Besides,
$\pi^\omega_*\omega$ is compatible with the cluster algebra transformations.
Moreover, if $\gamma$ is any other 2-form on the image $\Cal X^\omega$
compatible with the cluster algebra transformations, then its pullback
is a compatible 2-form on the preimage $\Cal X$. Thus, we get the following
result.

\proclaim{Corollary~2.2} Let $Z$ be irreducible, then there exists a
unique, up to a constant factor, $2$-form on $\Cal X^\omega$ compatible
with $\A(Z)$. Moreover, this form is symplectic.
\endproclaim

The unique symplectic form described in the above corollary is called
the {\it Weil-Petersson\/} form associated with the cluster algebra
$\A(Z)$. The reasons for this name are explained in the Introduction.

In our previous paper~\cite{GSV} we have studied Poisson structures
compatible with a cluster algebra. Fix a multiindex $I$ and consider the
reduced cluster algebra $\A_I(Z)$, as explained at the end of Section~2.1. By
Theorem~1.4 of~\cite{GSV}, there exists a linear space of Poisson structures
compatible with $\A_I(Z)$. To distinguish a unique, up to a constant factor,
Poisson structure, we consider the restriction of the form $\omega$ to the
affine subspace $f_i=\text{const}$, $i\notin I$. This restriction is
a symplectic form, since the kernels of $\omega$ are transversal to the above
subspace, and it is compatible with $\A_I(Z)\subset \A(Z)$. Its dual is
therefore a Poisson structure on this affine subspace, compatible with
$\A_I(Z)$. By the same transversality argument, we can identify each affine
subspace $f_i=\text{const}$, $i\notin I$, with the image of the quotient
map $\pi_Z$, and hence each of the obtained Poisson structures can be
identified with the dual to the Weil-Petersson symplectic form.

In ~\cite{GSV} we have considered an alternative description of the Poisson
structure dual to
$\pi^\omega_*\omega$. It follows from the proof of Theorem~1.4 in~\cite{GSV}
that it is the unique (up to a nonzero
scalar factor) Poisson structure $P$ on $\Cal X^\omega$ satisfying the
following condition: for any Poisson structure on $\Cal X$
compatible with any reduced cluster algebra $\A_I(Z)$, the projection
$p\:\Cal X\to \Cal X^\omega$ is Poisson. The basis $\tau_i$
is log-canonical with respect to $P$,
and $\{\tau_i,\tau_j\}=\lambda z_{ij}\tau_i\tau_j$
for some constant $\lambda$.

\heading{3. Main example: Teichm\"uller space}
\endheading

Let $\Sigma$ be a Riemann surface of genus $g$
with $s$ punctures (marked points).

\subheading{3.1. Teichm\"uller space, horocycles, decorated
Teichm\"uller space} Let $\TE_g^s$ denote the Teichm\"uller space
of marked complete metrics on $\Sigma$ having constant negative
curvature $-1$ and a finite area, modulo the action of the
connected component of identity in the group of 
diffeomorphisms of $\Sigma$. By the uniformization theorem, $\Sigma$
is identified with the quotient of the complex upper half-plane
modulo the action of a discrete M\" obius group. A horocycle
centered at a point $p$ at the absolute is a circle orthogonal to
any geodesic passing through $p$. All horocycles centered at some
fixed point $p$ can be parametrized by a positive real number
called the {\it height\/} of the horocycle. The {\it decorated\/} Teichm\"uller
space $\pTE_g^s$ classifies hyperbolic metrics on $\Sigma$ with a
chosen horocycle around each of the marked points; $\pTE_g^s$ is a
fiber bundle over $\TE_g^s$ whose fiber is $\R_+^s$.

More exactly, let $V$ be a
Minkowsky $3$-space, i.e. a real vector space of dimension~$3$ with a 
nondegenerate quadratic form $\langle\cdot,\cdot\rangle$ of type $(2,1)$.
Choosing appropriate coordinates, one can write the Minkowsky metric as
$-dv_0^2+dv_1^2+dv_2^2$. Consider the two-sheeted hyperboloid 
$\{v\in V\: \langle v,v\rangle=-1\}$, and denote by $\H$ its upper sheet.   
It is well known that $\H$ is isometric to the Poincar\'e disk model for the
hyperbolic plane. An explicit isometry is given by the radial projection 
from the origin to the unit disk $D$ about $(1,0,0)$. For $x,y\in\H$,
the Poincar\'e distance $d$ between the projections of $x$ and $y$ to
$D$ satisfies $\cosh^2d=\langle x,y\rangle^2$.

The light cone $L$ is defined by $L=\{v\in V\: \langle v,v \rangle =0\}=
\{v\in V\: v_0^2=v_1^2+v_2^2\}$. The positive light cone $L^+$ is
given by $L^+=\{v\in L\: v_0>0\}$. The radial projection extends to a map
$\pi\:\H\cup L^+\to D\cup S^1_\infty$, where $S^1_\infty$ is the boundary 
circle of $D$. A point $w\in L^+$ corresponds to the
horocycle $\{x\in\H\:\langle w,x\rangle=-1\}$ of height $w_0$ about $\pi(w)$.

Let us recall here a construction of coordinates on the decorated
Teichm\"uller space (see~\cite{Pe1}). 

Fix any ideal geodesic
triangulation $\Delta$ with vertices at the marked points. Recall
that an ideal geodesic triangulation 
is a maximal family of disjointly embedded simple geodesic arcs 
subject to the condition that no complementary region of $\Delta$
in $\Sigma$ is a monogon or a bigon (see Fig.~1). An ideal triangulation
subdivides $\Sigma$ into open
pieces (triangles) which are images of ideal triangles in the
complex upper-half plane under the uniformization map. Observe that
$\Delta$ may contain triangles with two coinciding sides (see Fig.~1). 
Any edge $e\in\ \Edge(\Delta)$ of
this triangulation has an infinite hyperbolic length, however the
hyperbolic length of its segment between the corresponding
horocycles is finite. Let us denote by $l(e)$ the signed length of
this segment, i.e., we take the length of $e$ with the positive
sign if the horocycles do not intersect, and with the negative
sign if they do, and define $f(e)=\exp(l(e)/2)$. Equivalently, a choice 
of a geodesic between two marked points corresponding to an edge $e$
together with two horocycles 
determines two points $x,y$ on the positive light cone, and 
we put $f(e)=\sqrt{-\langle x,y\rangle}$.

\vskip 15pt
\centerline{\hbox{\epsfxsize=10cm\epsfbox{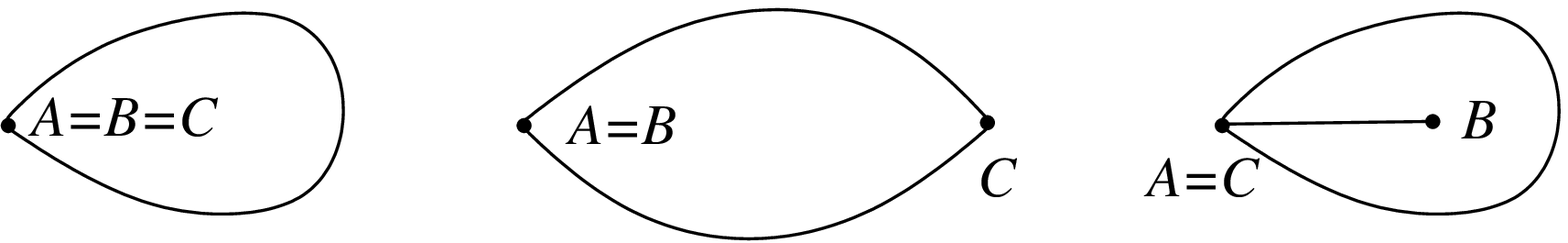}}}
\midspace{0.1mm} \caption{Fig.~1. Monogon, bigon, and a triangle with two
coinciding sides }

The following result is well known (see~\cite{Pe1}).

\proclaim{Theorem 3.1} For any ideal triangulation $\Delta$, the
functions $f(e)$, $e\in\ \Edge(\Delta)$, define a homeomorphism
between the decorated Teichm\"uller space $\pTE_g^s$ and
$\R^{6g-6+3s}_+$.
\endproclaim

The coordinates described in Theorem~3.1 are called the {\it
Penner coordinates\/} on~$\pTE_g^s$.

\subheading{3.2. Whitehead moves and the cluster algebra
structure} Given a triangulation of $\Sigma$, one can obtain a new
triangulation via a sequence of simple transformations called
flips, or Whitehead moves, see Fig.~2. Note that some of the sides
$a, b, c, d$ can coincide, while $p$ is distinct from any of them.

\vskip 15pt
\centerline{\hbox{\epsfysize=3,5cm\epsfbox{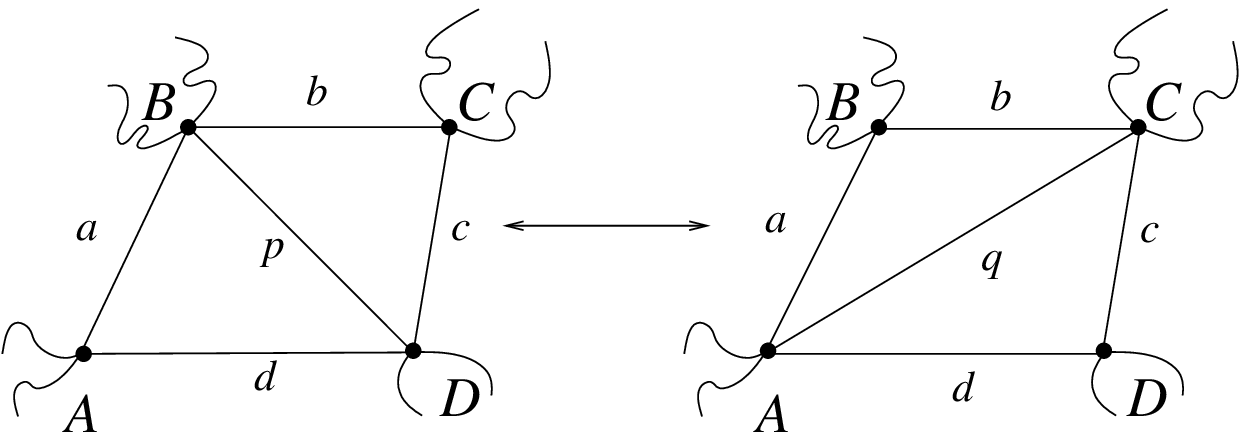}}}
\midspace{0.1mm} \caption{Fig.~2. Whitehead move }

\noindent Under such a transformation the Penner coordinates
undergo the following change described by the famous {\it Ptolemy
relation\/}:
$$
  f(p)f(q)=f(a)f(c)+f(b)f(d). \tag 3.1
$$

Note that this transformation looks very much like a cluster
algebra transformation; we would like to make this statement more
precise. Let us call an ideal triangulation of $\Sigma$ {\it
nice\/} if all the edges of any triangle are pairwise distinct. We
call a nice triangulation {\it perfect\/} if, moreover, all
vertices have at least three incident half-edges (in other words,
loops are counted with multiplicity $2$). We will consider Penner
coordinate systems only for nice triangulations. Note that every
edge in a perfect triangulation borders two distinct triangles
with different edges. So, $a$ and $c$ can coincide with each
other, but not with any  of $b$ and $d$; similarly, $b$ and $d$
can coincide as well.

 Therefore in this situation the Ptolemy
relation holds for any edge. In general, flips do not preserve the
perfectness of a triangulation; however, the result of a flip of a
perfect triangulation is a nice one.

Given a perfect triangulation $\Delta$, we construct the following
coefficient-free cluster algebra $\A(\Delta)$: the Penner
coordinates are cluster coordinates, and the transformation rules
are defined by Ptolemy relations. The transformation matrix
$Z(\Delta)$ is determined by $\Delta$ in the following way.

Let $\nu_P(a,b)$ be the number of occurrences of the edge $b$
immediately after $a$ in the counterclockwise order around vertex
$P$. For any pair of edges $a,b\in \Edge(\Delta)$, put
$$
Z(\Delta)_{ab}=\sum_{P\in \Vert(\Delta)} (\nu_P(a,b)-\nu_P(b,a)),
$$
and define $\A(\Delta)=\A(Z(\Delta))$.

For example, the left part of Fig.~3 shows a triangulation of the sphere
with three marked points. Here $\nu_P(a,b)=\nu_P(b,a)=1$ and
$\nu_Q(a,b)=\nu_Q(b,a)=\nu_R(a,b)=\nu_R(b,a)=0$, hence
$Z(\delta)_{ab}=0$. The right part of Fig.~3 shows a triangulation of the
torus with one marked point. Here $\nu_P(a,b)=0$, $\nu_P(b,a)=2$, and
hence $Z(\Delta)_{ab}=-2$.

\vskip 15pt
\centerline{\hbox{\epsfysize=3cm\epsfbox{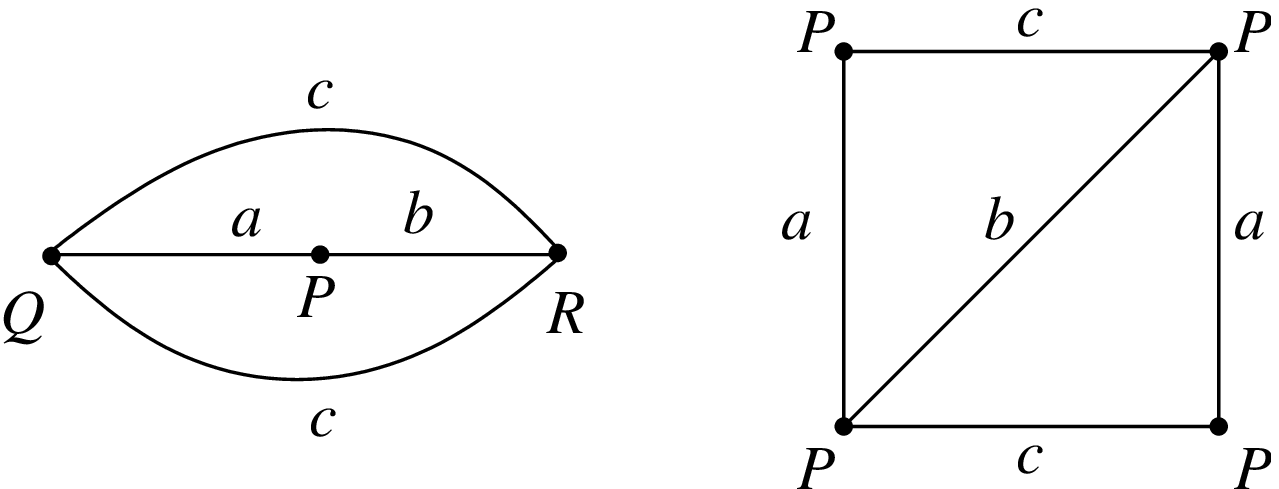}}}
\midspace{0.1mm} \caption{Fig.~3. To the definition of $Z(\Delta)$ }

 Let $\Delta\mapsto \bar\Delta$ be a flip shown on Fig.~2.
By the definition of $\A(\Delta)$ we see that cluster variables in
$\A(\Delta)$
corresponding to the cluster obtained from the initial one by the
transformation $T\: f(p)\mapsto f(q)$ coincide with Penner coordinates
for the flipped triangulation $\bar\Delta$.
To check that the flip is indeed a cluster algebra
transformation, it is enough to check that the transformation matrix
$\bar Z$
is determined by the adjacency of edges in $\bar\Delta$ in the same way as
$Z=Z(\Delta)$ above.

To prove this,  recall that the new transformation matrix $\bar Z$
is obtained from $Z$ by the mutation rule~(2.2).
Assume first that the degrees of both endpoints of $p$ are at least three;
the Whitehead move in this case is shown
on Fig.~2.

Note that the image of the diagonal $p$ under the flip is the diagonal $q$.
Using formula~(2.2) we obtain
$\bar z_{ad}=z_{ad}-1$, $\bar z_{ba}=z_{ba}+1$, $\bar z_{cb}=z_{cb}-1$,
$\bar z_{dc}=z_{dc}+1$, $\bar z_{xq}=-z_{xp}$, where $x$ runs over
$\{a,b,c,d\}$.

For instance,
substituting the values $z_{da}=z_{ap}= z_{pd}=z_{cp}=z_{pb}=z_{bc}=1$,
$z_{ab}=z_{cd}=0$ into (2.2) we obtain
$\bar z_{ab}=\bar z_{bq}= \bar z_{qa}=\bar z_{dq}=\bar z_{qc}=\bar z_{cd}=1$
and $\bar z_{ad}=\bar z_{bc}=0$ which corresponds to the adjacency rules
of Fig.~2. Similarly, one can easily check that the adjacency rules hold for
other triangulations containing the rectangle $ABCD$.

This shows that Penner coordinates for two perfect triangulations
related by a Whitehead move
form two sets of cluster coordinates for adjacent clusters in $\A(\Delta)$.

A problem arises in the case of non-perfect triangulations, as
shown in Fig.~4.

\vskip 15pt
\centerline{\hbox{\epsfxsize=8cm\epsfbox{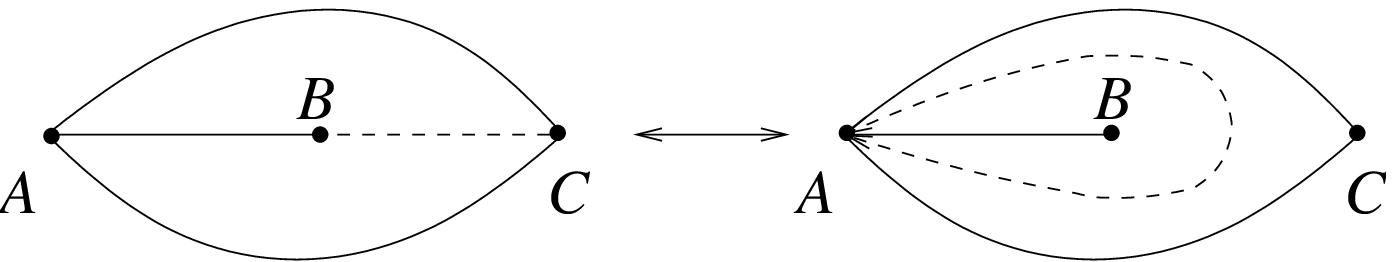}}}
\midspace{0.1mm} \caption{Fig.~4. Prohibited Whitehead move }

Indeed, Ptolemy relation implies $f_{BC} f_{A\bigcirc
A}=f_{AB}(f_{A\frown C}+f_{A\smile C})$ However, by the definition
of cluster algebras, the right part of the cluster algebra
relation must be irreducible, and in fact the cluster algebra
relation are $f_{BC} f_{A\bigcirc A}=f_{A\frown C}+f_{A\smile C}$.
We thus see that the flip with respect to the edge $BC$ breaks the
cluster algebra rules, meaning that Penner coordinates after the
flip differ from the coordinates prescribed by the cluster algebra
rules.

To overcome this problem we are going to prohibit flips of certain
diagonals in non-perfect triangulations. Namely, we call a
Whitehead move {\it allowed\/} only if both endpoints of the
flipped diagonal before the move have degrees at least three. As
we have already mentioned, the transformation rule for Penner
coordinates coincides with the transformation rule for the
corresponding cluster algebra, and we are guaranteed that Penner
coordinates of triangulations obtained by a sequence of such flips
form clusters in the same cluster algebra $\A(\Delta)$.

Observe that all nice triangulations might be obtained one from
another by a sequence of allowed flips. This can be proved by an
immediate modification of the proofs provided in~\cite{Ha, Bu},
where the authors consider the same property for ideal
triangulations~\cite{Ha}, or all triangulations~\cite{Bu}; another
proof follows from the fact that nice triangulations label maximal
simplices of triangulations of the Teichm\"uller space,
see~\cite{Iv} and references therein. Therefore, cluster algebras
$\A(\Delta)$ defined by all perfect triangulations $\Delta$ of a
fixed punctured surface are isomorphic to the same
cluster algebra, which we denote by $\A(\Sigma)$. 

The clusters of $\A(\Sigma)$ corresponding to nice triangulations
are called {\it Teichm\"uller clusters}. Let $\X(\Sigma)$ be the
cluster manifold for the cluster algebra $\A(\Sigma)$. It is easy to see
that the decorated
Tecim\"uller space $\tilde \TE_g^s$ coincides with the positive part of 
any toric chart defined by a Teichm\"uller cluster.

\subheading{3.3. The Weil-Petersson form}
Recall that if $\Delta$ is a perfect triangulation, then flips of all edges
are allowed. Let us define the {\it star\/} of $\Delta$ as the subset of
clusters in $\A(\Sigma)$ formed by $\Delta$ itself and all the clusters
obtained by flips with respect to all edges of $\Delta$; we
denote the star of $\Delta$ by $\star(\Delta)$. 
It follows immediately from
the proof of Theorem~2.1 and the connectedness of the graph of $\Delta$
that there exists a unique (up to a constant) differential 2-form on 
$\X(\Sigma)$ compatible with
all clusters in $\star(\Delta)$; moreover, this form is compatible with the
whole algebra $\A(\Sigma)$. Recall that
this form is given in coordinates $f(e)$, $e\in\Edge(\Delta)$, by the
following expression:
$$
\omega=\text{const}\cdot\sum_{b\curvearrowleft a}
\frac{df(a)\wedge df(b)}{f(a)f(b)}, \tag 3.2
$$
where $b\curvearrowleft a$ means that edge $b$ follows immediately after
$a$ in the counterclockwise order around some vertex of $\Delta$.

It is well-known that the Teichm\"uller space $\TE_g^s$ is a symplectic
manifold with respect to the famous Weil-Petersson symplectic form $W$.
Following ~\cite{Pe2}, we recall  here the definition of $W$.
The cotangent space $T^*_\Sigma\TE_g^s$ is formed by holomorphic quadratic
differentials $\phi(z)dz^2$. Define the Weil-Petersson nondegenerate co-metric
by the following formula:
$$
\langle \phi_1(z) dz^2,\phi_2(z) dz^2\rangle =
\frac{i}{2}\int_\Sigma \frac{\phi_1(z)\phi_2(z)}
{\lambda(z)}dz\wedge\overline{dz},
$$
where $\lambda(z)|dz^2|$ is the hyperbolic metric on $\Sigma$. The
Weil-Petersson metric is defined by duality. Moreover, the Weil-Petersson
metric is K\"ahler. Hence, its imaginary part is a symplectic
$2$-form $W$, which is called the Weil-Petersson symplectic form.

As we have already mentioned above, the decorated Teichm\"uller space is 
fibered
over $\TE_g^s$ with a trivial
fiber $\R_{>0}^s$. The projection $\rho\:\tilde \TE_g^s\to \TE_g^s$
is given by forgetting the horocycles.

Given a geodesic triangulation, we introduce {\it Thurston shear 
coordinates\/} on the
Teichm\"uller space. Following~\cite{Th, Pe1, Fo}, we associate to
each edge of the triangulation a real number. Choose an edge $e$ and
two triangles incident to it and lift the resulting rectangle to
the upper half plane. Vertices of this geodesic quadrilateral lie
on the real axis. We obtain therefore four points on the real axis,
or more precisely, on $\R P^1$. Among the constructed four points
there are two distinguished ones, which are the endpoints of the edge $e$
we have started with. Using the M\"obius group action on the upper
half plane, we can shift these points to zero and infinity,
respectively, and one of the remaining points to $-1$. Finally, we assign to
$e$ the logarithm of the coordinate of the fourth point (the
fourth coordinate itself is a suitable cross-ratio of those four points).

To show that we indeed obtained coordinates on the Teichm\"uller space,
it is
enough to reconstruct the surface. We will glue the surface out of
ideal hyperbolic triangles. The lengths of the sides of
ideal triangles are infinite, and therefore we can glue two
triangles in many ways which differ by shifting one triangle
w.r.t. another one along the side. The ways of gluing triangles can be
parametrized by the cross-ratio of four vertices of the obtained
quadrilateral (considered as points of $\R P^1$).

By~\cite{Pe1, Fo}, the projection $\rho$ is written in these
natural coordinates as follows:
$$
g(e)=\log f(e_1)+\log f(e_3)-\log f(e_2)-\log f(e_4),
$$
where $\{g(e)\: e\in \Edge(\Delta)\}$ and $\{f(e)\: e\in
\Edge(\Delta)\}$ are Thurston shear coordinates on the Teichm\"uller space
and Penner coordinates on the decorated
Teichm\"uller space, respectively, with respect to the same
triangulation $\Delta$; notation $e_i$ is explained in
Fig.~5.

\vskip 15pt
\centerline{\hbox{\epsfysize=3cm\epsfbox{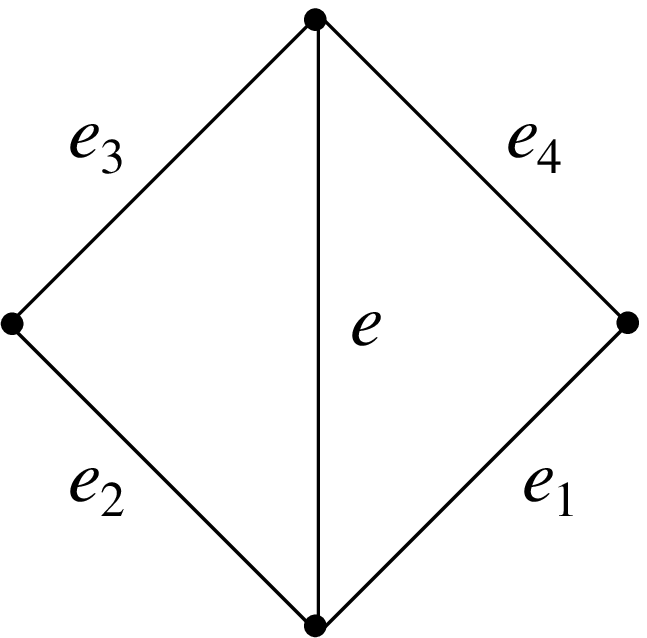}}}
\midspace{0.1mm} \caption{Fig.~5. Notation for neighboring edges}

We see that $\exp(\rho)$ coincides with the transition
from cluster variables to $\tau$-coordinates in the cluster
algebra $\A(\Sigma)$.
Therefore, the pullback of the Weil-Petersson symplectic form
determines a degenerate $2$-form of corank $s$ on the decorated
Teichm\"uller space. For any nice triangulation $\Delta$ its expression
in terms of Penner coordinates $f(e)$, $e\in\Edge(\Delta)$, is as
follows:
$$
\omega=\frac12\sum_{e\in\Edge(\Delta)}\left(dx_{e_1}\wedge dx_e -
dx_{e_2}\wedge dx_e+dx_{e_3}\wedge dx_e - dx_{e_4}\wedge dx_e\right),
$$
where $x_e=\log f(e)$, see for instance~\cite{FR}.

Comparing this expression with (3.2) we can summarize the previous
discussion as follows.

\proclaim{Theorem 3.2} {\rm(i)} The decorated
Teichm\"uller space $\tilde \TE_g^s$ is the subset of the cluster manifold
$\X(\Sigma)$ given by the positive part of the
toric chart defined by any Teichm\"uller cluster.
 
{\rm(ii)} The projection to the $\tau$-coordinates is the
forgetful  map from the decorated Teichm\"uller space to the
Teichm\"uller space.

{\rm(iii)} The corank of the corresponding transformation
matrix equals $s$. 

{\rm(iv)} The unique compatible symplectic structure
induced on the Teichm\"uller space is the classical Weil-Petersson
symplectic structure.
\endproclaim

In the next section, using the cluster algebra approach, we
describe subsystems of Penner coordinates such that restriction of
$\omega$ to these subsystems is nondegenerate. In particular, we
obtain an independent combinatorial proof of Theorem~3.2(iii).
Choosing complementary coordinates as stable variables and
following the procedure suggested in~\cite{GSV}, one obtains the
Weil-Petersson Poisson structure on the decorated Teichm\"uller
space.

\subheading{3.4. Gauge group action}
Given a triangulation $\Delta$ of $\Sigma$, consider the transformation
matrix $Z(\Delta)$ for the cluster algebra $\A(\Delta)$. Our goal is to prove
that the corank of $Z(\Delta)$ equals $s$
and to choose a subset
of Penner coordinates such that the restriction of $Z$ onto these coordinates
is nondegenerate.

Let us prove first that $\corank Z(\Delta)\gs s$.
Indeed, the gauge group $G$ of horocycle scalings acts
on the decorated Teichm\"uller space. Namely, we can multiply the height
of the horocycle about point $P$ by a positive factor $\lambda(P)$.
We thus get an action of the  multiplicative group $G=(\R_+)^s$  on
the decorated Teichm\"uller space, which is described in Penner
coordinates as follows:
$$
\left(\lambda_{P_1},\dots,\lambda_{P_s}\right)\:f(e)\mapsto
\lambda_{P_i}\cdot\lambda_{P_j}\cdot f(e), \tag 3.3
$$
where $e=P_iP_j\in\Edge(\Delta)$, see~\cite{Pe1}.

Consider an arbitrary cluster algebra $\A(Z)$.
Assume that a $n$-tuple of integer weights $w_u=(w_{u,1},\dots,w_{u,n})$ is
given at any
vertex $u$ of the tree $\Bbb T_n$. We define
a {\it local toric action\/} on the cluster $C_u$ at $u$ as the
$\R^\star$-action given by the
formula $\{f_{u,1},\dots,f_{u,n}\}\mapsto \{f_{u,1}\cdot
t^{w_{u,1}},\dots,
f_{u,n}\cdot t^{w_{u,n}}\}$. We say that local toric actions are
{\it compatible\/} (in $\A(Z)$) 
if for any two clusters $C_a$ and $C_b$ the following diagram is
commutative:
$$
\CD
C_a @>>> C_b\\
@V t^{w_a}VV @V t^{w_b}VV\\
C_a @>>> C_b
\endCD
$$
where the horizontal arrows are defined by cluster transformations
(2.1). In this case, local toric actions together define a
global toric action on $\A(Z)$. In other words,
local toric actions extend to a global one if all
relations (2.1) are homogeneous with respect to this action.

It is easy to see that local toric actions (3.3) are compatible in 
$\A(\Sigma)$, and hence they can be extended to a global toric action.
Lemma~2.3 in~\cite{GSV} claims that the dimension of the group of all
global toric actions
is equal to the corank of the transformation matrix $Z$. Thus
$\corank Z(\Delta)\gs s$.
Let us prove that $\corank Z(\Delta)=s$.

Let $\Delta$ be a nice triangulation of a genus $g$ surface $\Sigma$ with
$s$ punctures. Let $V(\Delta)$ be the vector space generated by the basis
$\{x_e\: e\in\Edge(\Delta)\}$, where $x_e=\log f(e)$,
$G$ be the gauge group acting on $V(\Delta)$ by~(3.3).
The form $\omega$ can be written, up to a constant, as
$\sum_{b\curvearrowleft a} d x_a\wedge d x_b$.
 For any subset $S$ of $\Edge(\Delta)$ we put
$V(S)=\spa\{x_e\:e\in S\}\subseteq V(\Delta)$. We call a subset
$S\subseteq\Edge(\Delta)$
{\it representative\/} if the gauge group $G$ acts on $V(S)$ faithfully
and $S$ is minimal with respect to this property. We denote by $\Gamma(S)$ 
the spanning subgraph of $\Delta$  with the edge set $S$.

\proclaim{Theorem 3.3} {\rm (i)}
Let $\Delta$ be a nice triangulation of $\Sigma$,
$S$ be any representative subset of $\Edge(\Delta)$, $R$ be the
complement of
$S$ in $\Edge(\Delta)$. If $R\ne\emptyset$, the restriction of the form
$\omega$ onto $V(R)$ is a nondegenerate $2$-form.

{\rm (ii)}
Each connected component of $\Gamma(S)$ contains
exactly one cycle and possibly a number of trees attached to it.
The length of this cycle is odd.
\endproclaim

\demo{Proof} (i) Let $V_0(\Delta)$ denote the vector space generated by
the basis $\{y_P\: P\in\Vert(\Delta)\}$, and define a map
$\partial^*\: V_0(\Delta)\to V(\Delta)$ by $y_P\mapsto \sum_{e=PQ}x_e$.
Evidently, $\partial^*$ is an injection. We want to prove that the complex
$V_0(\Delta)@>\partial^*>> V(\Delta)@>\zeta^*>> V(\Delta)$ is exact, 
where the second map is given by the matrix 
$Z^*(\Delta)=-Z(\Delta)$.
Equivalently, it is enough to prove that the dual complex 
$V(\Delta)@>\zeta>> V(\Delta) @>\partial>> V_0(\Delta)$
is exact. 

The map $\partial$ is given by $x_e\mapsto y_P+y_Q$
for $e=PQ$, so $\ker\partial$ is generated by the vectors $\sum_{e\in
\Edge(\Delta)}\alpha_ex_e$
such that $\sum_{e=PQ}\alpha_e=0$ for any point $P\in\Vert(\Delta)$. 
Consider the image of the map $\zeta$ given by 
$Z(\Delta)$. For any 
$e\in\Edge(\Delta)$ define the {\it alternating quadrangle\/} 
$q(e)\in V(\Delta)$ by $q(e)=Z(\Delta)x_e$; evidently, the image of $\zeta$
is spanned by $\{q(e)\: e\in\Edge(\Delta)\}$.

Let $C=(e_1,\dots,e_{2k})$ be an even
length cycle in the graph of $\Delta$ with a specified first edge $e_1$.
Define $\alt(C)\in V(\Delta)$ to be $\sum(-1)^ix_{e_i}$. Clearly, 
$\partial(\alt(C))=0$. We claim that both $\ker\partial$ and $\Im\zeta$ are 
spanned by $\alt(C)$, where $C$ runs over all even length cycles in the
graph of $\Delta$. 

Let us show first that any $\alt(C)$ is a linear combination of alternating 
quadrangles. We start from the following particular case.
Let $PQR$ and $PMN$ be two triangles of $\Delta$ sharing the vertex $P$.
Define $b_P(Q,R;M,N)\in V(\Delta)$ by $b_P(Q,R;M,N)=\alt(PQRPMN)=
x_{PQ}+x_{PR}-x_{QR}-
x_{PM}-x_{PN}+x_{MN}$. To prove that $b_P(Q,R;M,N)$ is a sum of alternating 
quadrangles, we proceed as follows. Order the edges incident to $P$ cyclically
in such a way that $e_1=PQ$, $e_2=PR$, $\dots$, $e_{l-1}=PM$, $e_l=PN$. 
It is easy
to see that $b_P(Q,R;M,N)=\sum_{i=2}^{l-1}q(e_i)$.

For an arbitrary even length cycle $C=(e_1,\dots,e_{2k})$ we proceed as 
follows. Let $e_i=P_iP_{i+1}$; pick $Q_i$ in such a way that $P_iP_{i+1}Q_i$
is a triangle of $\Delta$. Then $2\alt(C)=\sum(-1)^ib_{P_i}(P_{i-1},Q_{i-1};
P_{i+1},Q_i)$, and we are done.

Further, we fix a vector  $x=\sum_{e\in
\Edge(\Delta)}\alpha_ex_e$
such that 
$$
\sum_{e=PQ}\alpha_e=0\qquad \text{for any $P\in\Vert(\Delta)$}, \tag3.4
$$
and prove that it can be represented
as a linear combination of $\alt(C)$ for even length cycles $C$. Consider
the subgraph $\Gamma$ of the graph of $\Delta$ containing only edges with 
$\alpha_e\ne0$.
Assume first that $\Gamma$ contains a simple even length cycle $C$, that is, 
all the vertices of $C$ are distinct. Define 
$\alpha=\min_{e\in C}|\alpha_e|$, then either $x+\alpha \alt(C)$ or
$x-\alpha \alt(C)$ satisfies condition (3.4), while the corresponding
graph has at least one edge less than $\Gamma$. If $\Gamma$ does not contain a
simple even length cycle, then $\Gamma$ is a tree of edges and simple odd 
cycles,
that is, any two simple odd cycles have at most one common vertex
(see \cite{We, Ex.~4.2.18}). By the faithfulness of the action, 
$\Gamma$ itself is not a simple odd cycle.
Moreover, condition (3.4) implies that $\Gamma$ does not contain
pendant vertices, hence it contains at least two simple odd cycles.
Evidently, for any two simple odd cycles $C_1$ and $C_2$ in $\Gamma$ 
there exist
vertices $P_1\in C_1$ and $P_2\in C_2$ and a simple path $\pi$ with the
endpoints $P_1$ and $P_2$ that does not have common edges with $C_1$ and $C_2$
(in particular, it may occur that $P_1=P_2$ and $\pi$ is an empty path).
We then define $\alpha_i=\min_{e\in C_i}|\alpha_e|$ for $i=1,2$,
$\alpha_\pi=\frac12\min_{e\in \pi}|\alpha_e|$, set $\alpha=\min\{\alpha_1,
\alpha_2,\alpha_\pi\}$ and proceed as before with the non-simple even length 
cycle $C$ consisting of $C_1$, $C_2$, and $\pi$, the latter traversed in both 
directions. As a result, we again obtain a vector satisfying (3.4) such
that the corresponding graph has at least one edge less than $\Gamma$, and the
statement follows by induction.

It remains to prove that $V(R)$ intersects $\ker\zeta^*$ trivially. Indeed,
let $x\in \ker\zeta^*$. By above, this means that 
$x=\partial^*(\sum_{P\in\Vert(\Delta)}\lambda_Py_P)$. If, in addition, 
$x\in V(R)$, then $\{\lambda_P\:P\in\Vert(\Delta)\}$ acts trivially on
$V(S)$, which contradicts to the faithfulness of the action.

(ii) Follows immediately from the faithfulness of the action.
\qed
\enddemo

\subheading{3.5. Denominators of transformation functions and
intersection numbers}
In conclusion, we would like to prove the following result 
that gives a nice interpretation of denominators of
transition functions in terms of intersection points on the surface
$\Sigma$. Let us express all cluster variables in $\A(\Delta)$ as
rational functions in the cluster variables of the initial cluster. 
Recall that all the clusters in $\A(\Delta)$ are related to nice
triangulations of $\Sigma$, while cluster variables correspond to
edges of this triangulation.  Abusing notation, we denote an edge of a
triangulation and the corresponding cluster variable by the same letter.
For a cluster variable $p$ and an initial cluster variable $x$ we
denote by $\delta_x(p)$ the exponent of $x$ in the denominator of $p$. 
Besides, let $a$ and $b$ be two edges possibly belonging to different 
nice triangulations; we denote by $[a,b]$ the number of their inner 
intersection points (that is, the ones distinct from their ends).

\proclaim{Theorem 3.4}
$\delta_x(p)=[x,p]$.
\endproclaim

\demo{Proof} Consider an arbitrary edge $p$ and two triangles bordering on
$p$ (see Fig.~6). Denote by $q$ the edge obtained by flipping $p$. We want to
prove the following two relations:
$$\align
[p,x]+[q,x]&=\max\{[a,x]+[c,x],[b,x]+[d,x]\},\tag 3.5\\
\delta_x(p)+\delta_x(q)&=\max\{\delta_x(a)+\delta_x(c),
\delta_x(b)+\delta_x(d)\}.
\tag 3.6\endalign
$$

To prove (3.5), we consider all the inner intersection 
points of $x$ with the sides $a$, $b$, $c$, $d$. These points break $x$
into consecutive segments, each either lying entirely inside the quadrangle
$abcd$, or entirely outside this quadrangle. We distinguish three types of
segments lying inside the quadrangle. The segments of the first type
intersect two adjacent sides of the quadrangle and exactly one of the
diagonals $p$ and $q$ (see Fig.~6). Thus, the contribution of any such 
segment to both sides of (3.5) equals $1$. The segments of the second type
intersect two opposite sides of the quadrangle and both diagonals. Moreover,
all segments of the second type intersect the same pair of the opposite sides,
otherwise $x$ would have a self-intersection. This means that each segment
of the second type contributes $2$ to both sides of~(3.5). Finally, the
segments of the third type start at a vertex and intersect the opposite
diagonal and one of the opposite sides of the quadrangle. They are treated 
in the same way as the segments of the second type, the only difference being 
that their contribution to both sides of~(3.5) equals $1$. To get (3.5)
is enough to sum up the contributions of all the segments.

\vskip 15pt
\centerline{\hbox{\epsfxsize=3cm\epsfbox{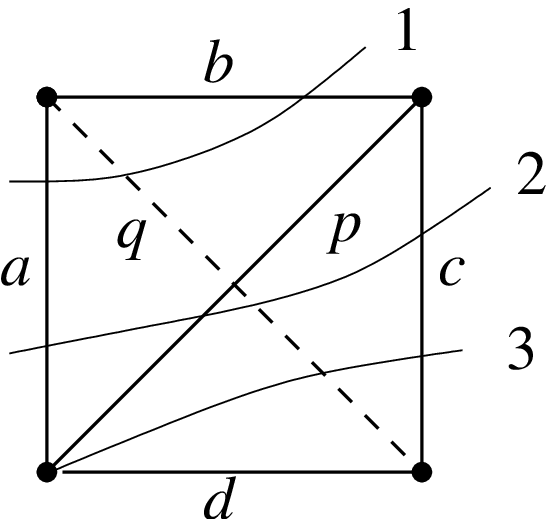}}}
\midspace{0.1mm} \caption{Fig.~6. Types of segments}

To prove (3.6), consider an arbitrary cluster variable, say $p$, as a function
of the initial cluster variable $x$, provided all the other initial cluster
variables are fixed to~$1$. Then $p=P(x)/x^k$, where $P(x)$ is a polynomial 
in $x$ with a nonzero constant term, and $k\gs 0$. It is easy to notice 
that the above constant term is positive. Indeed, by Theorem~3.1, $p>0$
for any choice of positive values of the initial cluster coordinates. On
the other hand, the sign of $p$ coincides with the sign of the constant term
under consideration, provided $x$ is sufficiently small. Therefore,
$\delta_x(ac+bd)=\max\{\delta_xa+\delta_xc, \delta_xb+\delta_xd\}$, and
(3.6) follows immediately from this fact and the exchange relation 
$pq=ac+bd$. 

Thus, we have proved that both $\delta_x(p)$ and $[p,x]$ satisfy the same
relation with respect to flips. To get the statement of the theorem 
it suffices to notice that it holds evidently when $p$ and $e$ belong to
two adjacent clusters.
\qed
\enddemo

This observation was independently made by Dylan Thurston, who
used it, in particular, to describe relations between Penner
coordinates corresponding to distant triangulations.

We illustrate the lemma with the following example, see Fig.~7.

\vskip 15pt
\centerline{\hbox{\epsfxsize=10cm\epsfbox{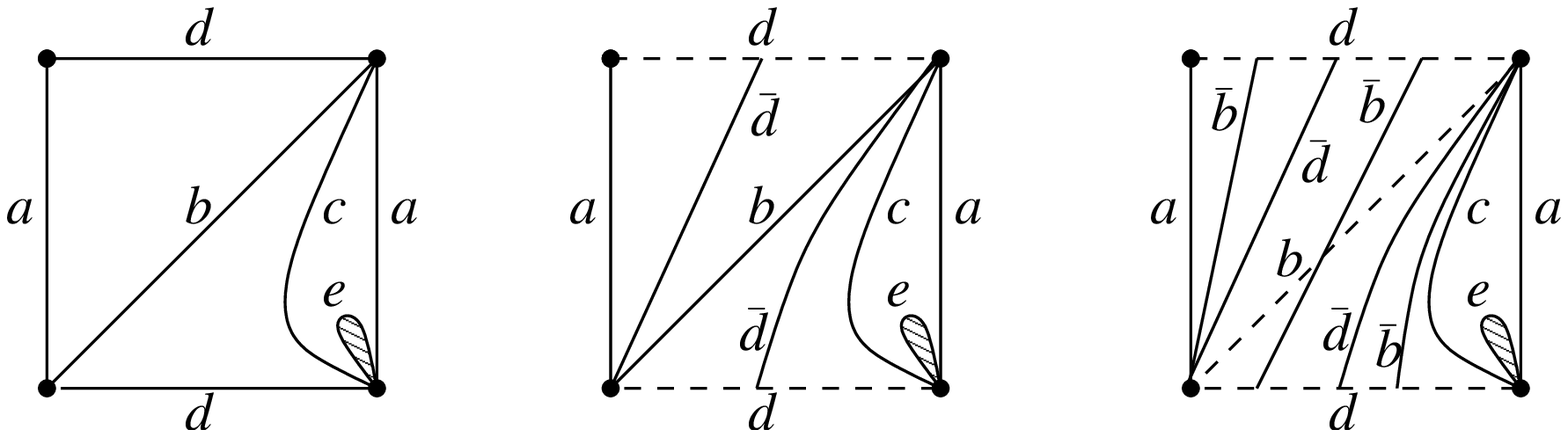}}}
\midspace{0.1mm} \caption{Fig.~7. Denominators and the number
of intersections }

After the first move we have $\bar d=(ac+b^2)/d$, and hence
$\delta_d(\bar d)=1$,
which corresponds to the intersection of the edges $d$ and $\bar d$. After
the second move we have $\bar b=(acd^2+(ac+b^2)^2)/bd^2$ and hence
$\delta_b(\bar b)=1$, $\delta_d(\bar b)=2$, which corresponds to the
intersection
of the edge $\bar b$ with $b$ and to two intersections of $\bar b$ with $d$.

\end{document}